\input AHTOH-E.STY
\baselineskip 10.9pt

\UDC{
512.624.3 
+ 512.552   
}


\title{%
Balanced factorisations
in some algebras
}
\author{%
Anton A. Klyachko
\quad
Andrey M. Mazhuga
\quad
Anastasia N. Ponfilenko
}
\address{
\myAddress
\quad
mazhuga.andrew@yandex.ru
\quad
stponfilenko@gmail.com
}

\grantsFirstSecond{\RFBR15-01-05823}

\abstract{%
We prove that, in any field of characteristic not two and not three except 
$\F_5$, each element decomposes into a product of four factors whose sum 
vanishes. We also find all $k,n,q$ such that every $n\times n$ matrix over 
the $q$-element field decomposes into a product of $k$ commuting matrices 
whose sum vanishes.}

\s 0.
Introduction

Suppose that an element $a$ of a ring is factored into a product:  
$a=a_1a_2\dots a_k$; we call this factorisation \emph{balanced} if 
$\sum a_i=0$. 
It is easy to show (see [KV16]) that, 
\disp{\sl in any 
field of characteristic not two, each element admits a balanced 
decomposition into a product of $k$ factors for each $k\ge5$. \(While, for 
each $k<5$, this assertion is no longer valid.\)} 
We prove the following theorem answering thereby two question from [KV16] 
(one of which had been known earlier~[Iva13]).

\Th 1.
In any field of characteristic not two and not three except the 
five-element field $\F_5$, each element admits a balanced decomposition 
into a product of four factors; if the field is infinite, then each 
element admits infinitely many such decompositions.

\noindent
We are unsure whether this result was known earlier, because
\-
A high-school teacher of the second author, Dmitrii Vitalievich Andreev, 
in 2003 suggested such a problem on classes (for the field 
of rationals); the problem was not solved, but the teacher's
hints make the second author suspect (now, when he already knows 
a solution) that Dmitrii Vitalievich was able to solve this problem; 
\- 
in 2016 (already knowing a solution), the second author asked such a
question (again for the field of rationals) on the forum 
{\tt http://math.stackexchange.com/} and shortly someone posts 
a solution, different from ours, but probably also valid; now this 
question is, unfortunately, removed and we cannot give an exact 
references.

\enditem
The following theorem from [KV16] describes, for each $k$, all finite 
fields, in which each element admits a balanced decomposition into 
a product of $k$ factors.

\Th on balanced factorisations in finite fields {\rm [KV16]}.
Suppose that $k\ge2$ is an integer and $F$ is a finite field.
Then
any element of $F$ has a balanced factorisation  
into a product of $k$ factors if and only if
\itemitem{either}
$|F|=2$ and $k$ is even,
\itemitem{or}
$|F|=4$ and $k\ne3$,
\itemitem{or}
$|F|$ is a power of two but neither two nor four (and $k$ is arbitrary),
\itemitem{or}
$|F|\in\{3,5\}$ and $k\notin\{2,4\}$,
\itemitem{or}
$|F|=7$ and $k\notin\{2,3\}$,
\itemitem{or}
$|F|$ is neither a power of two nor three nor five nor seven and
$k\ne2$.
\enditem
{\rm In other words, the situation in finite fields is the following:}


{
\tabskip1em plus 2em
\halign to \hsize{
#\hfil&\strut\hfil#\hfil&\hfil#\hfil&\hfil#\hfil&\hfil#\hfil&\hfil#\hfil
\cr
          &$k=2$&$k=3$&$k=4$&$k=5,7,9,\dots$&$k=6,8,10,\dots$
\cr
$\F_2$    & yes & no  & yes &     no        &  yes
\cr
$\F_3$    & no  & yes & no  &     yes       &  yes
\cr
$\F_4$    & yes & no  & yes &     yes       &  yes
\cr
$\F_5$    & no  & yes & no  &     yes       &  yes
\cr
$\F_7$    & no  & no  & yes &     yes       &  yes
\cr
$\F_8,
\F_{16},
\F_{32},
\F_{64},
\dots$    & yes & yes & yes &     yes       &  yes
\cr
$\F_9,
\F_{11},
\F_{13},
\F_{17},
\dots$    & no  & yes & yes &     yes       &  yes
\cr
}
}

\smallskip

\medskip

\centerline{Table 1}

\medskip

The following theorem complements this result from [KV16].

\Th2.
{\let\itemitem\
Let $k,n\ge2$ be integers and $F$ be a finite field.
Any $n\times n$ matrix over $F$ decomposes into a
product of $k$ commuting matrices whose sum vanishes
if and only if
\newline
\itemitem{either}
$k=3$ and $|F|=5$,
\itemitem{or}
$k=3$ and $|F| \geqslant 8$,
\itemitem{or}
$k=4$ and $|F|=4$,
\itemitem{or}
$k=4$ and $|F| \geqslant 7$,
\itemitem{or}
$k \geqslant 5$ and $|F| \geqslant 3$.
}
\rm
In other words, the situation in the matrix algebras over finite fields is 
the following (ignore the superscripts for now):

\medskip\goodbreak

{
\tabskip1em plus 2em
\halign to \hsize{
#\hfil&\strut#\hfil&#\hfil&#\hfil&#\hfil&#\hfil
\cr
          &$k=2$    &$k=3$    &$k=4$        &$k=5,6,7,8,\dots$
\cr
$\F_2$    & no$^1$ & no$^0$ & no$^2$     &no$^{2}$
\cr
$\F_3$    & no$^1$ & no$^8$ & no$^0$     &yes$^{3,4}$
\cr
$\F_4,
\F_7$     & no$^1$ & no$^0$ & yes$^{5,7}$  &yes$^{3,4,5,7}$
\cr
$\F_5$    & no$^1$ & yes$^5$  & no$^0$     &yes$^{3,4}$
\cr
$\F_8,
\F_9,
\F_{11}
\F_{13}
\F_{16},
\dots$    & no$^1$ &yes$^{5,6}$& yes$^{5,7}$  &yes$^{3,4,5,7}$
\cr
\cr
}
}

\smallskip
\centerline{Table 2}

\noindent
Note that the size of matrices does not affect the answer, provided this
size is at least two.

In Section 1, we prove Theorem 1
and also another theorem on finite fields (Theorem 3), which plays the key 
role in the proof of Theorem 2.  
All argument in Section 1 is elementary, except
that the proof of Theorem 3 is substantially based on results
of [KV16] (that in turn are based on the theory of elliptic
curves). In Section 2, we prove Theorem 2.

{\smallskip\noindent\bf Notation},
which we use are standard. Note only that the symbol
$\F_q$ denotes the  $q$-element field and the letter $E$ always denotes
the identity matrix.

\s 1.
Fields

\noindent{\bf Proof of Theorem 1.} 
Look at this identity
$$
x=
\frac{2(1-4x)^2}{3(1+8x)}
\cdot
\frac{-(1+8x)}{6}
\cdot
\frac{-(1+8x)}{2(1-4x)}
\cdot
\frac{18x}{(1-4x)(1+8x)}.
$$
We can only suggest readers to verify 
this identity and the fact that the sum of factors is zero. In
the exceptional cases, when denominators vanishes, i.e. for
${x\in\{{1\over4},-{1\over8}\}}$, 
we can multiply $x$ 
by~$y^4$ choosing $y$ such that 
$xy^4\notin\{{1\over4},-{1\over8},0\}$ (this is possible in any field of 
characteristic not two and not three except $\F_5$), write a similar 
decomposition for $xy^4$ and then divide each factor by $y$:  
$$ 
x= 
\frac{2(1-4xy^4)^2}{3y(1+8xy^4)}
\cdot
\frac{-(1+8xy^4)}{6y}
\cdot
\frac{-(1+8xy^4)}{2y(1-4xy^4)}
\cdot
\frac{18xy^4}{y(1-4xy^4)(1+8xy^4)}.
$$
Thus we obtain a balanced decomposition of any element
into a product of four factors.
Moreover, the last identity 
implies that there are infinitely many such decompositions if the field is 
infinite, because of the following 
elementary fact (whose proof is left for readers as an easy exercise):  
\disp{\sl 
any non-constant rational fraction over an infinite field takes 
infinitely many values.} 
(Note that, for each $x$, the second factor in the last 
identity is a non-constant rational fraction 
$f(y)$.) 
This 
completes the proof. The five-element field is indeed an
exception, see Table~1.

\medskip

Similarly, we can obtain infinitely many
balanced decompositions into products of any larger number
of factors, i.e.
\disp{\sl in any infinite field characteristic not two each element
admits infinitely many balanced decompositions into products of $k$
factors for each $k\ge5$.}
For example, the following identity is obtained by a slight modification
of an identity from [KV16]; this gives infinitely many
balanced decompositions of any nonzero element of an infinite
field of characteristic not two into products of 2017 factors:
$$
x={xy^{2016}\over2}\cdot
{xy^{2016}\over2}\cdot
\(-xy^{2016}\)\cdot
{2\over xy^{2018}}\cdot
\(-{2\over xy^{2018}}\)\cdot
\({1\over y}\)^{1006}\cdot
\(-{1\over y}\)^{1006}.
$$

\bigskip

A decomposition  $x=x_1x_2\dots x_k$ is called \emph{power}
if all factors are equal: ${x_1=\dots=x_k}$ [KV16].

\Th 3.
{\let\itemitem\
Let $k\ge2$ be an integer and let $F$ be a finite field.
Any element of $F$ admits a non-power
balanced
decomposition into a product of $k$ factors
if and only if
\newline
\itemitem{either}
$k=3$ and $|F|=5$,
\itemitem{or}
$k=3$ and $|F| \geqslant 8$,
\itemitem{or}
$k=4$ and $|F|=4$,
\itemitem{or}
$k=4$ and $|F| \geqslant 7$,
\itemitem{or}
$k \geqslant 5$ and $|F| \geqslant 3$.
\rm In other words, the answer here is the same as in Theorem 2 (Table 2).
}

\Proof
Superscripts of the particles \emph{yes} and \emph{no} in Table 2 indicate
references to the cases below.

\Case 0:
no, because in these cases for some elements there are no
balanced decompositions (by Theorem on balanced factorisations
in finite fields).

\Case 1:
$k=2$~--- no.
For any element $a$, consider its balanced decomposition $a=xy$,
$x+y=0$. If the characteristic is 2, then the decomposition $a=(-x)x$
is power; if the characteristic is not two, then
not any element is square and, therefore, not any element
admits a balanced decomposition into a product of two factors.

\Case 2: $|F|=2$~--- no.
The decomposition of 1 can comprises only 1s and, therefore, it is
power.

\Case 3:
$k=5+2n$, where $n\ge0$ and $\Char F\ne2$ --- yes.
Let us use the universal formula for balanced decompositions
into a product of $5+2n$ factors from [KV16]:
$$
\pm a=(-a)\cdot{\frac a2}\cdot{\frac a2}\cdot{\frac 2a}\cdot{\frac {-2}a}
\cdot1^n\cdot(-1)^n
\qbox{(for $a\ne0$).}
$$
$\Char F\ne2$, therefore, $\frac2a\ne-\frac2a$
and, hence, this decomposition is non-power.
The zero element has an obvious balanced non-power decomposition:
$0=(-1)\cdot1\cdot0^{k-2}$.

\Case 4:
$k=6+2n$, where $n\ge0$ and $\Char F\ne2$ --- yes.
Let us use a formula for balanced decompositions
into a product of $6+2n$ factors from [KV16]:
Consider $c\in F$ such that $c^2\ne a$
(such $c$ exists, except for the case where $F=\F_3$
and $b=1$; but in this case everything is obvious).
Put
$b=\frac{c^2-a}c$. Then
$$
\pm a=(-c)\cdot
(c-b)
\cdot
{\frac b 2}
\cdot
{\frac b 2}
\cdot
{\frac 2 b}
\cdot
{\frac{-2} b}\cdot
1^n \cdot (-1)^n .
$$
Since $\Char F\neq2$, we have $\frac2b\ne-\frac2b$ and,
therefore, the decomposition is non-power.

\Case 5:
In these cases, balanced decompositions exist by Theorem on
balanced factorisations in finite fields; these decompositions cannot  
be power, because the number of factors is not divisible by the 
characteristic.

\Case 6:
$|F|\ge9$, $\Char F\ne2$, and $k=3$ --- yes.

Consider two cases.
If $\Char F\ne3$, then the Theorem
on
balanced factorisations in finite fields
provides us with
a balanced decomposition for any element;
since $\Char F\ne3$, this decomposition cannot be power.

To prove the assertion for characteristic three, we need a lemma.

\Lemma.
In the field $\F_{3^n}$, where $n\ge2$, there exists a nonzero square
$u$ such that $u+1$ is also a nonzero square.

\Proof
If 2 is a square, then $u=1$ is the required element. Otherwise,
suppose that, for any $u\notin\{0,1,2\}$, 
the following holds:
\disp{
if $u$ is a square, then $u+1$ is not a square.}
Then, any
set of the form $\{u,u+1,u+2\}$ contains at most one square.
Therefore, the number of squares is at most $3^{n-1}+1$. On the other 
hand, in the field of characteristic 3, there are precisely
$\frac{3^n+1}2$ squares. This implies the inequality 
${\frac{3^n+1}2\le3^{n-1}+1}$ that holds only for $n=1$. This 
contradiction completes the proof of the lemma.

\medskip

Let us resume the proof of Theorem 3.
We have to show that, in a finite field of characteristic three
and cardinality at least nine, each
element has a balanced non-power decomposition into a product of three 
factors.

Note that, in such a field, any element is a cube.
We want to find a balanced decomposition of an element
$a=b^3\ne0$:
$$
a=xyz,
\quad
x+y+z=0.
\qqbox{Eliminating $z$, we obtain}
yx^2+y^2x+a=0.
\eqno{(*)}
$$
Let us solve this equation
with respect to $x$.
By lemma, there exists $\tau^2\ne0$ such that $\tau^2+1\ne0$ is
also a square: $\tau^2+1=\pi^2\ne0$.
Take $y=\frac{b+b\pi}{2}$.
Note that $y\ne b$, because the equality $y=b$ would mean that
$\pi=1$ and $\tau=0$.
Therefore, the discriminant of the quadratic equation $(*)$ is a square:
$$
D=y^4-4ay=y(y^3-b^3)=y(y-b)^3=
\frac{b\pi+b}{2}\cdot\(\frac{b\pi-b}{2}\)^3=
(b^2(1+\tau^2)-b^2)(b\pi-b)^2=
b^2\tau^2(b\pi- b)^2
$$
and
equation $(*)$
has a solution. The obtained decomposition is not 
power, because $y^3\ne a$ (since $y\ne b$).

It remains to find a non-power balanced decomposition of zero, but this 
is an easy task:  $0 = (-1)\cdot1\cdot0$.

\Case 7:
$|F|=2^n$ and $k=4+2m$, where $m\ge0$.
Since $\Char F=2$, any element is a square:
$a=b\cdot b$.
If $a\ne1$, then $b\ne1$ and $a=b^2\cdot1^{2m+2}$ is a required
decomposition.
If $a=1$, then
$a=1=c^2\cdot\(\frac1c\)^2\cdot1^{2m}$ is a required
decomposition, where $c$ is any element different from 0 and 1.

\Case 8: $|F|=3$ and $k=3$ --- no.
The decomposition of 1 cannot can contain 0 and, therefore, must contain 1 
and $-1$, to make decomposition non-power. Then, the third 
factor must be zero, because the decomposition is balanced. 
This contradiction completes the proof of Theorem 3.

\s 2.
Matrices

{\noindent\bf Proof of Theorem 2.}
First note that, for $k=2$, Theorem 2 is valid:
\disp{\it the Jordan cell $J$ with eigenvalue zero and size
$n\times n$ is not a square in the ring of $n\times n$ matrices if 
$n\ge2$} 
(we leave the proof of this fact to readers as an exercise) and, 
therefore, the matrix $-J$ has no balanced decompositions into a product 
of two factors.

In the case $k\ge3$, Theorem 3 implies that it suffices to prove the 
following statement.

\Th 2$'$.
Let $n\ge2$ and $k\ge3$ are integers and let $F$ be a field (not 
necessarily finite). Then the following conditions are equivalent:  
\item{a)}
each matrix $n\times n$ over  $F$ has a balanced decomposition into a
product of $k$ commuting factors;
\item{b)}
each element of $F$ has a non-power balanced decomposition into a
product of $k$ factors.

\Proof

{\noindent\bf The implication $b)\imp a)$} 
follows immediately from the 
following fact proven in [KV16]:  
\disp{\it Let $F$ be a field and let $k$ be a
positive integer larger than two. If, in all finite extensions of 
the field~$F$, each element has a non-power balanced decomposition into a  
product $k$ elements, then the same is true for each element of each 
finite-dimensional associative algebra with unit over~$F$.}

{\noindent\bf The implication $a)\imp b)$.}
First, note that $0\in F$ has a balanced non-power decomposition 
into a product of $k$ factors for any $k\ge3$:
\quad
$0=0^{k-2}\cdot1\cdot(-1)$.
To obtain a decomposition of a nonzero element $a\in F$,
we need the following simple fact from linear algebra (the proof
is left to readers as an exercise):
\disp{\it The centraliser of the nilpotent Jordan cell $J$ of size
$n\times n$ in the algebra of ${n\times n}$ matrices consists of 
polynomials in $J$, i.e.  
$C(J)=\{a_0E+a_1J+\dots+a_{n-1}J^{n-1}\;|\;a_i\in F\}.$} 
Thus, if the Jordan cell $aE+J$ has a balanced decomposition 
$aE+J=X_1\dots X_k$ into a product of commuting matrices, then all 
matrices $X_i$ lie in the centraliser of $J$ and, by virtue of 
the fact mentioned above, we obtain a balanced decomposition of 
$a$ in $F$:  
$$ 
a=x_1\dots x_k, 
\qbox{where $x_i$ is the (unique) 
eigenvalue of $X_i$.} 
$$ 
It remains to note that this decomposition cannot 
be power for $a\ne0$. Indeed, assuming the contrary, we would obtain a 
balanced decomposition in the ring of matrices:  
$$ 
aE+J=(xE+J_1)\dots(xE+J_k), 
\qbox{ where $J_i$ are nilpotent 
commuting matrices.} 
$$ 
The balancedness of this decomposition means that 
$k$ is divisible by $\Char F$ and $\sum J_i=0$. But then, 
multiplying out brackets, we obtain 
$$ 
aE+J=(xE+J_1)\dots(xE+J_k)=aE+f(J_1,\dots,J_k), 
$$ 
where the polynomial $f$ 
has no terms degree 1, i.e.  the right-hand side of this equality 
is a matrix of the form $aE+J'$, where $(J')^{n-1}=0$.  This  
contradiction completes the proof of Theorems $2'$ and 2.

\bigskip

Generally (without the commutativity condition), the question on balanced
decompositions of matrices over finite fields remains open.

\Question.
For which $q$, $k$, and $n$, it is true that any $n\times n$ matrix 
over the $q$-element field has a balanced decomposition into a product
of $k$ matrices?

We can say only the following.

\item{1.}
In some cases the decompositions exist by Theorem 2.

\item{2.}
For $k=2\le n$, the decomposition does not exist, e.g., because 
factors such balanced decomposition must commute.

\medskip
\enditem
Moreover,
computer experiments show the following facts.

\item{3.}
Over the two-element field, 
all $2\times2$ matrices, except
$\pmatrix{
1 & 1 \cr
1 & 0 \cr
}$
and
(a similar matrix)
$\pmatrix{
0 & 1 \cr
1 & 1 \cr
}$,
admit balanced decompositions into a product of three factors, while
these two matrices have no such decompositions.

\item{4.}
The matrix
$\pmatrix{
1 & 0 \cr
1 & 1 \cr
}$
and two similar matrices over the two-element field have no
balanced decompositions into a product of four factors, while all other
$2\times2$ matrices over $\F_2$ have such decompositions.

\item{5.}
The matrix
$\pmatrix{
1 & 0 & 0 \cr
1 & 1 & 0\cr
0 & 0 & 1\cr
}$
and similar matrices over the two-element field have no
balanced decompositions into a product of three factors, while the 
remaining $3\times3$ matrices over $\F_2$ have such decompositions.

\item{6.}
All $3\times3$  matrices over $\F_2$ have balanced decompositions into
products of four factors.

\item{7.}
All $2\times2$ matrices over $\F_3$, $\F_4$, $\F_5$, and $\F_7$ have
balanced decompositions into products of three and four factors.
This implies that the same is true for any larger number of factors, 
because we can increase the number of factors by two multiplying 
decompositions by $E$ and $-E$.


\REFERENCES



\[Iva13]
Ivanishchuk A. V.
The experience of learning and research activity of students
in Lyceum 1511 (MEPhI) (in~Russian)
//
published
in the book
{\it Sgibnev A. I.
Research problems for beginners.
Moscow: MCCME,~2013.}
(Freely available at {\tt http://www.mccme.ru/free-books/})

\[KV16]
Klyachko A. A., Vassiliyev A. N.
Balanced factorisations //
arXiv:1506.01571.






\end